\documentclass[12pt,leqno]{article}

\usepackage{amsmath,amsfonts,amssymb,mathtools,amsthm,epsfig}
\usepackage{mathrsfs}
\usepackage{comment}

\usepackage{indentfirst}
\usepackage{hyperref}
\hypersetup{colorlinks=true,urlcolor=blue,linkcolor=blue,citecolor=blue}

\usepackage{enumitem}

\usepackage{psfrag}
\usepackage[mathscr]{euscript}
\usepackage[utf8]{inputenc}
\usepackage[english]{babel}
\usepackage[title]{appendix}

\usepackage{epsfig}

\usepackage{url}
\usepackage{setspace}

\usepackage{verbatim}
\usepackage[all]{xy}

\theoremstyle{plain}

\newtheorem{global-theorem}{Theorem}
\newtheorem{theorem}{Theorem}[section]
\newtheorem{lemma}[theorem]{Lemma}

\newtheorem{corollary}[theorem]{Corollary}

\theoremstyle{definition}

\newtheorem{remark}[theorem]{Remark}

\DeclareFontFamily{U}{rsf}{}
\DeclareFontShape{U}{rsf}{m}{n}{
  <5> <6> rsfs5 <7> <8> <9> rsfs7 <10->  rsfs10}{}
\DeclareMathAlphabet{\mathscr}{U}{rsf}{m}{n}
\newcommand{\mycal}[1]{\mathscr{#1}}

\DeclarePairedDelimiter\floor{\lfloor}{\rfloor}
\DeclarePairedDelimiter\fra{\{}{\}}

\newcommand{\op}[1]{\operatorname{#1}}

\newcommand{\bbC}{\mathbb{C}}
\newcommand{\bbP}{\mathbb{P}}
\newcommand{\bbQ}{\mathbb{Q}}
\newcommand{\bbZ}{\mathbb{Z}}
\newcommand{\cA}{\mathcal{A}}
\newcommand{\cO}{\mathcal{O}}
\newcommand{\net}{\mathbf{A}}
\newcommand{\discr}{\mathfrak{d}}
\newcommand{\bdelta}{\boldsymbol{\delta}}
\newcommand{\hi}{\boldsymbol{\imath}}
\newcommand{\muk}{\boldsymbol{v}}
\newcommand{\iu}{\sqrt{-1}}
\newcommand{\pt}{\mathsf{pt}}
\newcommand{\hp}{\mathsf{h}}
\newcommand{\Hp}{\mathsf{H}}

\newcommand{\Knum}{K_{\op{num}}}
\newcommand{\nc}{{\sf\bfseries nc}}
\newcommand{\bLambda}{\boldsymbol{\Lambda}}
\newcommand{\blambda}{\boldsymbol{\lambda}}

\DeclareFontFamily{U}{cyr}{}
\DeclareFontShape{U}{cyr}{m}{n}{
  <5> wncyr5 <6> wncyr6 <7> wncyr7 <8> wncyr8 <9> wncyr9 <10->
wncyr10}{}
\DeclareMathAlphabet{\mathcyr}{U}{cyr}{m}{n}
\input cyracc.def
\input epsf

\newcommand{\TSh}{{\mathcyr{\cyracc Sh}}}

\title{The intersection of three quadrics  in \texorpdfstring{$\bbP^{7}$}{P7} revisited}
\author{R. Donagi and T. Pantev}
\date{July 2024}

\begin{document}

\maketitle

\tableofcontents

\

\section{Introduction} \label{sec:intro}

In this paper we analyze certain additive invariants  of the sheaf theory of the natural gerbe associated with a very general  complete intersection of three quadrics in $\bbP^{7}$. Sheaf theories of $\mathbb{G}_{m}$-gerbes are very interesting for a variety of reasons, e.g. they contain essential information about the complexity of unramified Brauer classes, serve as a common playground for geometric constructions of representations and, most pertinently for us, provide  standard examples of non-commutative spaces (\nc-spaces for short).  

In the categorical approach to non-commutative  algebraic geometry over the complex numbers (see \cite{kkp,orlov-geometric.nc,maxim-ncgeometry} and references therein)  one studies cocomplete triangulated $\bbC$-linear dg categories. These categories are viewed as dg versions of the categories of quasi-coherent sheaves on the putative \nc-spaces and,  since in general  a direct independent description of such spaces is not available, they are taken as the primary geometric objects. The  categorical point of view towards \nc-geometry allows for unified treatment of complex and symplectic manifolds, as well as quantized and stacky geometries. Geometric properties, such as being algebraic, smooth,
proper, or Calabi-Yau, are notions intrinsic to \nc-spaces. Moreover the cohomology groups  of a smooth and proper \nc-space are enhanced with additional algebraic data, namely the data comprising a pure \nc-Hodge structure \cite{kkp}.

Classical algebraic geometry fits within this formalism, as a classical scheme or an algebraic stack $Y$ over $\bbC$ can be regarded as an \nc-space for which the corresponding category is $D_{\mathsf{qcoh}}(Y)$ -  a dg enhancement of the  quasi-coherent derived category of $Y$. There are many other examples of \nc-spaces coming from geometry or physics. A brief and incomplete list includes Fukaya categories of symplectic manifolds, matrix factorization categories, Fukaya-Seidel categories of a potential, and 
categories of boundary conditions of topological quantum field theories with $N=2$ susy. The properties and structures of an \nc-space corresponding to a cocomplete triangulated $\bbC$-linear category $\mathcal{C}$ are captured in the subcategory $\mathsf{Perf}(\mathcal{C})$ of compact objects in $\mathcal{C}$. For a classical quasi-separated scheme $Y/\bbC$ this notation is consistent with the standard notation, i.e. 
$\mathsf{Perf}(D_{\mathsf{qcoh}}(Y)) = \mathsf{Perf}(Y)$, where $\mathsf{Perf}(Y)$ is the full subcategory of objects locally quasi-isomorphic to finite complexes of locally free sheaves. 

An important class of \nc-spaces is comprised of the \nc-spaces of \emph{\bfseries geometric origin} \cite{orlov-geometric.nc}. By definition these are given  by the categories $\mathcal{C}$ for which  $\mathsf{Perf}(\mathcal{C})$ admits an embedding as an admissible subcategory in the category $\mathsf{Perf}(Y)$ for some scheme $Y/\bbC$. The \nc-spaces of geometric origin are particularly nice because  they are stable under natural noncommutative splitting and glueing constructions \cite{orlov-geometric.nc,alexncHPD}, and  because their  homological invariants are much more tractable than those of general \nc-spaces \cite{alexHC,alexDT}. Kuznetsov components of Fano manifolds, residual categories of spherical functors \cite{sashaCY,sasha.alex}, and derived categories of twisted sheaves on smooth projective varieties all provide important examples of \nc-spaces of geometric origin. 

\

These examples are of special  interest, since they are expected to give new birational invariants of smooth projective varieties \cite{ludmil-rationality,ikp,sasha-rationality}.
In particular, it is expected that truly non-commutative \nc-spaces of geometric origin will be birational invariants of the projective varieties through which they are realized geometrically. While general Kuznetsov components and residual categories are frequently  distinguished from smooth projective varieties by standard Hodge theoretic invariants, this is less common and harder when it comes to \nc-spaces associated to algebraic 
$\mathbb{G}_{m}$-gerbes, i.e. \nc-spaces corresponding to categories of twisted sheaves. 

To illustrate this point, recall that for any smooth projective variety $X/\bbC$, and any 
unramified Brauer class $\alpha \in \mathsf{Br}(X)$, we have a dg category $D^{b}(X,\alpha)$ of $\alpha$-twisted coherent sheaves, or equivalently of $\alpha$-twisted perfect complexes. This category embeds as an admissible subcategory in the Brauer-Severi variety of $\alpha$ and so corresponds to a \nc-space ${}_{\alpha}X$ of geometric origin. Unfortunately, it can happen that the \nc-space ${}_{\alpha}X$ is secretly commutative, even when the Brauer class $\alpha \in \mathsf{Br}(X)$ is non-trivial. Indeed, suppose that
$X$ is a smooth projective K3 equipped with a Weierstrass elliptic fibration. In this case one has a canonical identification $\mathsf{Br}(X) = \TSh(X)$ of the Brauer group of $X$ with the algebraic Tate-Shafarevich group  of the given elliptic fibration $\xymatrix@1{X
\ar[r]^-{\pi} & \bbP^{1} \ar@/^0.5pc/[l]^-{\sigma}}$. Furthermore, each element $\alpha \in \TSh(X)$ is realized by a smooth projective K3 $X_{\alpha}$ which is fibered by genus one curves and has Jacobian fibration isomorphic to $\xymatrix@1{X
\ar[r]^-{\pi} & \bbP^{1} \ar@/^0.5pc/[l]^-{\sigma}}$. The Poincar\'{e} sheaf on the fiber product ${}_{\alpha}X\times_{\bbP^{1}} X_{\alpha}$ induces an equivalence 
of categories $D^{b}(X,\alpha) \cong D^{b}(X_{\alpha})$ \cite{dp-torus.fibrations} and so the non-commutative space  ${}_{\alpha}X$ is isomorphic to the classical K3 $X_{\alpha}$. 
This phenomenon shows that the question of recognizing whether an \nc-space corresponding to a category of twisted sheaves is truly noncommutative  is quite delicate. 

A method for distinguishing such \nc-spaces from commutative ones was proposed 
by Kuznetsov in \cite{sashak-cubic}. It boils down to checking that the numerical $K$-lattice $(K_{\mathsf{num}}(Z),\chi_{Z})$ of a geometric \nc-space $Z$ can not be isomorphic to the numerical $K$-lattice of a smooth projective variety by checking that the lattice $K_{\mathsf{num}}(Z)$ cannot contain a class that behaves like the class of a point. In the same paper, Kuznetsov used this method to show that the category of twisted sheaves associated with a very general four dimensional cubic containing a plane is truly non-commutative.

In this paper we use Kuznetsov's method to show  that the natural  \nc-space attached to an intersection of three quadrics in $\bbP^{7}$ is truly non-commutative. More precisely, when $X$ is a general smooth complete intersection of three quadrics in $\bbP^{7}$, the derived category of $X$ has \cite{xie-quadric} a semiorthogonal decomposition 
\[
D^{b}(X) = \left\langle \, D^{b}(S,\alpha), \, \cO_{X}(1), \, \cO_{X}(2)  \, \right\rangle,
\]
where $S$ is a surface of general type, which is a double cover of the plane branched along a smooth octic and $\alpha$ is a $2$-torsion Brauer class on $S$. In order to analyze the properties of $D^{b}(S,\alpha)$, we consider 
the lattice 
\[
\bLambda := \left( \, H(S,B;\bbZ), \, \left\langle \, \sqrt{td(S)}\cdot(-), \, \sqrt{td(S)}\cdot(-) \right\rangle\, \right),
\]
where $B \in H^{2}(S,\frac{1}{2}\bbZ)$ is a class such that $\exp(2\pi\sqrt{-1}B) = \alpha \in H^{2}(S,\cO^{\times})$, $H(S,B;\bbZ) = \exp(B)(\oplus \, H^{p,p}(S,\bbQ)) \, \cap \, H^{\bullet}(S,\bbZ)$, and 
$\langle\, -\, , \, - \,\rangle$ is the Mukai pairing on $S$. Up to isomorphism the lattice $\bLambda$ depends only on $(S,\alpha)$ but not on $B$, and we check in Corollary~\ref{cor:obstruction} that $\bLambda$ can not be isomorphic to the Mukai lattice of integral Hodge classes of any smooth projective surface.
An immediate consequence is the following

\

\noindent
{\bfseries Theorem~\ref{theo:nc}} {\em \ Let $X$ 
be a very general intersection of $3$ quadrics in $\bbP^{7}$, and let 
$(S,\alpha)$ be the associated surface of general type with a $2$-torsion Brauer class. Then the associated numerical $K$-lattice $(\Knum(S,\alpha),\chi_{(S,\alpha)})$ is not isomorphic to $(\Knum(M),\chi_{M})$ for any smooth projective  surface $M$. In particular $D^{b}(S,\alpha)$ can not be equivalent to the derived category of any smooth projective surface.  
}

\

\noindent
The fact that $\bLambda$ is not isomorphic to the Mukai lattice of Hodge classes of a smooth projective surface 
has applications in birational geometry. 
The theory of Hodge atoms developed in \cite{kkpy-atoms,kkpy-enhatoms} implies that the $\alpha$-twisted  $\bbZ$-Hodge structure of $S$, viewed as a $\bbZ/2$-graded Hodge structure, together with its Mukai pairing gives an enhanced Hodge atom  $\blambda(X)$ in the atomic decomposition of $X$. If $X$ is very general, so that $\op{Pic}(S) = \bbZ\hp$, then one can check \cite{kkpy-enhatoms} that the  integral lattice of Hodge classes of the atom  $\blambda(X)$ is isomorphic to the lattice 
$\bLambda$. Since $\bLambda$ can not be the lattice of Hodge classes of any enhanced atom of a surface, the non-rationality criterion \cite[Proposition~5.21]{kkpy-atoms} implies the following

\

\noindent
{\bfseries Corollary.} {\em \ The very general intersection of three quadrics in $\bbP^{7}$ is not rational.}

\

\

\noindent
This corollary gives a new proof of the non-rationality result of Hassett-Pirutka-Tschinkel \cite{hpt}.

\

\noindent \textbf{Acknowledgments.} 
We would like to thank Ludmil Katzarkov, Maxim Kontsevich, Alex Perry, and Tony Yue Yu for valuable conversations on the subject of this note.  During the preparation of this work, Ron Donagi's research was supported in part by NSF grants DMS 2001673 and 2401422, by NSF FRG grant DMS 2244978, and by Simons HMS Collaboration grant 390287. 
Tony Pantev's research was supported in part by NSF/BSF grant DMS-2200914,  NSF FRG grant DMS 2244978, and by Simons HMS
Collaboration grant 347070.

\

\section{Geometric setup} \label{sec:setup}

Consider a smooth Fano fourfold $X$ which is a complete intersection of three quadrics in $\bbP^{7}$. In other words $X := \cap_{t \in \net} Q_{t}$ is the base locus of a general net of quadrics $Q_{t} \subset \bbP^{7}$, where $\net \cong \bbP^{2}$ denotes the parameter space of the net. We will assume that 
$\net \subset \left| \cO_{\bbP^{7}}(2) \right|$ is sufficiently general so that the  singular quadrics in this net all have corank one, and the discriminant curve $\discr \subset \net$ of the net will be a smooth curve of degree $8$.

\subsection{Birational fibration by quadric surfaces}  \label{ssec:qsurfaces}

\noindent
The complete intersection $X$ always contains lines and we fix a general line $\ell \subset X$ which does not intersect the curve in $\bbP^{7}$ which is the union of singularities of the singular quadrics of the net. 

Consider the incidence variety
\[
Y := \left\{ (t, \Pi)  \, \left| \,  \ell \subset \Pi \subset Q_t \, \right.\right\}  \subset  \net \times \mathsf{Gr}(3,8).
\]
This $Y$ fibers over $\net$. A general point $t \in \net$ determines first a $5$-dimensional projective subspace denoted $\mathsf{P}_t$:
\[
\mathsf{P}_t := \cap_{x \in \ell} T_{Q_{t},x}  \subset \bbP^{7},
\]
then a $3$ dimensional subquotient
\[
\mathsf{p}_t  :=  \mathsf{P}_t/\ell,
\]
and finally a quadric surface
\[
\mathsf{q}_t \subset \mathsf{p}_t,   \quad \mathsf{q}_t := (Q_t \cap \mathsf{P}_t)/\ell.
\]
The incidence variety  $Y$ fibers over $\net$, and the fiber over a general $t \in \net$ is isomorphic to the quadric surface $\mathsf{q}_t  \subset \mathsf{p}_t$. 

Furthermore, there is a birational  isomorphism $a: Y\to X$ from our $Y$ back to $X$. Starting with the pair $(t, \Pi)$, the intersection of $\Pi$ with a quadric $Q_s$ depends only on the image $\overline{s}$  of $s$ in the quotient line $\net/t$, and this intersection consists of a pair of lines: the fixed $\ell$ plus a variable line $\ell_{\overline{s}}$. The lines  $\left\{\ell_{\overline{s}}\right\}_{
\overline{s} \in \net/t}$ form a (linear) pencil of lines in $\Pi$, so they have a unique intersection point $a(t, \Pi)  \in \cap_{s \in \net} Q_{s} = X$. The inverse of $a$ 
is given by  $a^{-1}(x) = (t, \Pi)$, where $\Pi$ is the span of $x$ and $\ell$, and $t$ is the unique point such that $\Pi \subset Q_t$.

Note that our genericity assumptions on $\net$ and $\ell$ imply that $Y$ is smooth, and that the discriminat curve for the family 
$Y \to \net$ of quadric surfaces $\{\mathsf{q}_{t}\}_{t \in \net}$ is the same as the discriminant curve for the family 
$\mathcal{Q} \to A$ of six dimensional quadrics $\{Q_{t}\}_{t \in \net}$. In fact, if we start with $\ell \subset Q_{t}$ for a smooth quadric $Q_t$, then the resulting $\mathsf{q}_t$ is also smooth. Conversely, if $Q_t$ is singular and $\ell$ does not intersect its singular locus, then $\mathsf{q}_t$ is singular with the same singularity type. The family of two-dimensional quadrics $Y \to \net$ is the so-called hyperbolic reduction \cite{sashak-hyperbolic,xie-quadric} of the family of six-dimensional quadrics $\mathcal{Q} \to \net$ with respect to the isotropic family of lines $\ell \times \net \subset \mathcal{Q}$. Families related by hyperbolic reduction are known to share many invariants. In particular, it is known \cite{sashak-hyperbolic,xie-quadric} that their discriminant loci are isomorphic (as we saw in our example above) and that they have common pieces in the semiorthogonal decompositions in their derived categories.

For future reference, note that under our genericity assumptions, the total space $\mathcal{Q}$ of the family of six-dimensional quadrics over $\net$ is a divisor of bidegree $(1,2)$ in $A\times \bbP^7$. We can think of this as a family of $8\times 8$ symmetric matrices parametrized by $\net$, with entries linear forms on $\net$. The discriminant divisor $\discr\subset \net$ of the map $\mathcal{Q} \to \net$ is given by the determinant of the family,  which is a polynomial $\bdelta$ of degree $8$ on $\net$. 

To do the same calculation from the point of view of the two-dimensional quadrics $\mathsf{q}_{t}$, we have to first understand the family of three-dimensional projective spaces $\mathsf{p}_t$. This is the projectivization of a rank $4$ vector bundle $E$ on $\net$ and again our $Y$ is a smooth  divisor on $\mathbb{P}(E)$, which after appropriately normalizing $E$ by twisting it by a line bundle on $\net$, will be   in the linear system $p^{*}\mathcal{O}_{\net}(1) \otimes \mathcal{O}_{E}(2)$, where $\mathcal{O}_{E}(1)$ is the relative hyperplane bundle on $p : \mathbb{P}(E) \to \net$ that is characterized by the property $p_{*}\mathcal{O}_{E}(1) = E^{\vee}$. In these terms, the family of two-dimensional quadrics is a (symmetric) map $\sigma : E \to E^{\vee}\otimes \mathcal{O}_{\net}(1)$ and so the discriminant curve is a section in $\det(E)^{\otimes -2}\otimes \mathcal{O}_{\net}(4)$. In particular we must have 
$\det(E) = \cO_{\net}(-2)$. We can describe the bundle $E$ somewhat more concretely. For ease of reference, we record the following. 

\

\begin{lemma} \label{lem:ambient}
For a general  line $\mathbf{b} \subset \net$, we have 
$E|_{\mathbf{b}} \cong \cO_{\mathbf{b}}(-1)^{\oplus 2}\oplus \cO_{\mathbf{b}}^{\oplus 2}$.
\end{lemma}

\noindent
{\bfseries Proof:} \ Our $\bbP^{7}$ is a projectivization of an $8$-dimensional vector space $V$, and $\ell$ is the projectivization of  a $2$-dimensional subspace $L \subset V$. The family $\sigma$ of $8\times 8$ symmetric matrices (parametrized by $\net$, with entries linear forms on $\net$) all vanish on $L$. If we take the first two coordinates on $V$ to be the coordinates on $L$ (i.e. $L$ is given by vanishing of the remaining $6$ coordinates), we see that the upper left $2\times 2$ corner of these matrices vanishes identically, and then the lower left $2\times 6$ block gives the equations of the fibers of  $E$ inside the fixed vector space  $V/L$. In other words, $E$ fits in a short exact sequence:
\[
0 \to E \to \cO_{\net}^{\oplus 6} \to \cO_{\net}(1)^{\oplus 2} \to 0
\]
of vector bundles on $\net$.  In particular, the determinant of $E$ is $\cO_{\net}(-2)$. 

Now note that the  restriction of $E$  to a line $\mathbf{b} \subset \net$ will be isomorphic to 
$\oplus_{i = 1}^{4} \cO_{\mathbf{b}}(a_{i})$, and since $E|_{\mathbf{b}}$ is a subbundle in the trivial bundle $\cO_{\mathbf{b}}^{\oplus 6}$, we must have $a_{i} \leq 0$ and $\sum_{i=1}^{4} a_{i} = -2$. Thus 
$E|_{\mathbf{b}}$ has to be either $\cO_{\mathbf{b}}(-2) \oplus \cO_{\mathbf{b}}^{\oplus 3}$ or 
$\cO_{\mathbf{b}}(-1)^{\oplus 2} \oplus \cO_{\mathbf{b}}^{\oplus 2}$. These two are distinguished by the value of $h^0(\mathbf{b},E)$.
But $H^0(\mathbf{b},E)$ is the set of solutions of $2\times 2=4$ linear equations in $6$ variables, so for generic $\ell$ and $\mathbf{b}$  it is $2$ dimensional, and so the splitting type of $E$ is $E|_{\mathbf{b}} = \cO_{\mathbf{b}}(-1)^{\oplus 2} \oplus \cO_{\mathbf{b}}^{\oplus 2}$. \ \hfill \ $\Box$

\

\subsection{The \texorpdfstring{$2$}{2}-torsion gerbe} \label{ssec:2torsiongerbe}

\noindent
Consider next the double cover $\pi : S \to \net$ parametrizing the rulings on fibers of the quadric surface fibration $Y \to \net$. By construction  $S$ is  a smooth surface which is the double cover of $\net$ branched at the smooth octic curve $\discr \subset \net$, and the universal ruling of $Y \to \net$ is a \linebreak $\bbP^{1}$-bundle $\mycal{R} \to S$, i.e. a Brauer-Severi variety over $S$. Let $\alpha \in H^{2}(S,\cO^{\times})$ be the class of $\mycal{R}$ in the Brauer group $H^{2}(S,\cO^{\times})$ of $S$. By definition, $\alpha$ is a $2$-torsion element in  $H^{2}(S,\cO^{\times})$  which can be viewed also as the class of a sheaf $\cA \to S$ of rank $4$  Azumaya algebras over $S$. Note that for a very general choice of $\net \subset |\cO_{\bbP^{7}}(2)|$, the surface $S$ will be Noether-Lefschetz general and so we will have $\op{Pic}(S) = \bbZ$, generated by the hyperplane line bundle $\pi^{*}\cO_{\net}(1)$. We will write  
$\hp = \pi^{*}c_{1}(\cO_{\net}(1)) \in H^{1,1}(S)\cap H^{2}(S,\bbZ)$ for the corresponding hyperplane class.  Note that by construction this $\hp$ has self intersection  $\hp^{2} = 2$. Furthermore $S$  is a very ample divisor in the total space of $\cO_{\net}(4)$, and so, by the Lefschetz hyperplane section theorem for quasi-projective varieties \cite{TrangLefschetz,GM}, it will be simply connected. In particular, we have  $h^{1}(S,\bbZ) = h^{3}(S,\bbZ) = 0$. Therefore from the exponential sequence 
\[
\xymatrix@1@C+2pc@M+0.6pc{
0 \ar[r] & \bbZ \ar[r] & \cO_{S} \ar[r]^-{\exp(2\pi \iu(-))} & \cO^{\times} \ar[r] & 0
}
\]
we conclude that we can find  a class $\beta \in H^{2}(S,\cO_{S})$ such that $\exp(2\pi \iu \beta) = \alpha$. But $\alpha$ is $2$-trorsion, and so we can find a class $B \in H^{2}(S,\frac{1}{2}\bbZ)$, such that $\beta$ is the image of \linebreak $B \in H^{2}(S,\frac{1}{2}\bbZ) \subset H^{2}(S,\bbQ)$ under the map $H^{2}(S,\bbQ) \to H^{2}(S,\cO_{S})$ induced from $\bbQ \subset \cO_{S}$. 

Let us choose a class $B \in H^{2}(S,\frac{1}{2}\bbZ)$ such that $\exp(2\pi i B) = \alpha$. 
In general, the class $B$ itself is not unique: the different choices of $B$ are a coset under the subgroup $H^{2}(S,\bbZ) + H^{1,1}(S,\frac{1}{2}\bbZ) \subset H^{2}(S,\frac{1}{2}\bbZ)$.  Still, certain numerical invariants and numerical properties of $B$ are independent of the choice of an element in the coset and depend just on $\alpha$. Below we will explore these invariants in detail.

Let $\fra{t} := t - \floor{t}$ denote the fractional part of $t \in \bbQ$. Then it is easy to see \cite[Lemma~6.1]{sashak-cubic} that the fractional part $\fra{B\hp}$ does not depend on the choice of $B$. Indeed if $u \in H^{2}(S,\bbZ)$, then $(B+u)\hp = B\hp + u\hp$, and since $u\hp$ is an integer, we have 
$\fra{(B+u)\hp} = \fra{B\hp}$. Similarly $\fra{(B+\frac{1}{2}\hp)\hp} = \fra{B\hp + \frac{1}{2}\hp^{2}} = \fra{B\hp + 1} =\fra{B\hp}$. In fact $\fra{B\hp}$ can be computed in terms of the geometry of the $\bbP^{1}$-bundle $\mycal{R} \to S$ defining $\alpha$. The following is \cite[Lemma~6.2]{sashak-cubic}:

\

\begin{lemma} \label{lem:degree}
Let $C \subset S$ be a smooth curve in the linear system $|\hp|$. Then $\mycal{R}|_{C}$ is isomorphic to $\bbP(W)$ for some rank two vector bundle $W$ on $C$, and
\[
\fra{B\hp} = \left\{ \frac{1}{2} \deg W \right\},
\]
for any such $W$.
\end{lemma}

\

\noindent
From here we get an immediate 

\

\begin{corollary} \label{cor:{Bh}}
We have $\fra{B\hp} = \frac{1}{2}$ and $\fra{B^{2}}$ is independent of the choice of $B$.
\end{corollary}
\noindent {\bfseries Proof.} This is an adaptation of the analysis in \cite[Lemma~6.4]{sashak-cubic} to our setting.

Let $\mathbf{b} \subset \net$ be a line transversal to the degree $8$ curve $\discr \subset \net$. Then $C = \pi^{*}\mathbf{b}$ is a smooth hyperelliptic curve of genus $3$ in the linear system $|\hp|$. Abusing notation, we will still write $\pi : C \to \mathbf{b}$ for the covering map. By definition  $\mycal{R}|_{C}$ is the Hilbert scheme parametrizing lines in the quadrics $\mathsf{q}_{t}$, $t \in \mathbf{b}$. We will write $E_{\mathbf{b}}$, $Y_{\mathbf{b}}$ for the restrictions $E|_{\mathbf{b}}$, $Y|_{\mathbf{b}}$, etc. As before, the restricted family $Y_{\mathbf{b}}$ is a smooth $3$-fold fibered by quadric surfaces over the line $\mathbf{b}$. It is a smooth divisor on 
$\bbP(E_{\mathbf{b}})$ in the linear system $\left| \cO_{E_{\mathbf{b}}}(2)\otimes p^{*}\cO_{\mathbf{b}}(1)\right|$, given by a section 
$\sigma_{\mathbf{b}} \in H^{0}(\mathbf{b},\op{Sym}^{2}E_{\mathbf{b}}^{\vee} \otimes \cO_{\mathbf{b}}(1)) = H^{0}(\bbP(E_{\mathbf{b}}),\cO_{E_{\mathbf{b}}}(2)\otimes p^{*}\cO_{\mathbf{b}}(1))$.

By Lemma~\ref{lem:ambient} we have that $E_{\mathbf{b}} = \cO_{\mathbf{b}}^{\oplus 2}\oplus \cO_{\mathbf{b}}(-1)^{\oplus 2}$, and so we can consider 
the rank $3$ subbundle $F := \cO_{\mathbf{b}}^{\oplus 2} \oplus \cO_{\mathbf{b}}(-1) \subset E_{\mathbf{b}}$. Then $\bbP(F) \subset \bbP(E_{\mathbf{b}})$ is a family of hyperplanes in the fibers of $\bbP(E_{\mathbf{b}})$ and so intersects $Y_{\mathbf{b}}$ in a surface $\mathfrak{C}$ which fibers over $\mathbf{b}$ with fiber over $t \in \mathbf{b}$ a conic $\mathfrak{C}_{t} = \mathsf{q}_{t}\cap \bbP(F_{t}) \subset \mathsf{p}_{t} = \bbP(E_{t})$. For general $\net$, $\mathbf{b} \subset \net$, the surface $\mathfrak{C}$ will be smooth and all conics $\mathfrak{C}_{t}$ will be reduced. Furthermore, the conic bundle $\mathfrak{C} \to \mathbf{b}$ is given by the restricted family of quadratic equations $\sigma_{\mathbf{b}}|_{F} : F \to F^{\vee}\otimes \cO_{\mathbf{b}}(1)$ and so its discrimant is a section in 
$\det(F)^{\otimes (-2)} \otimes \cO_{\mathbf{b}}(3) = \cO_{\mathbf{b}}(5)$. In other words, the conic bundle $\mathfrak{C} \to \mathbf{b}$ has $5$ reducible fibers. Again, by taking $\mathbf{b}$ and $\net$ general, we can ensure that 
the discriminant of the family of conics $\mathfrak{C} \to \mathbf{b}$
is disjoint from the discriminant of the family of quadric surfaces $Y_{\mathbf{b}} \to \mathbf{b}$. Altogether we get $13 = 8 + 5$
special points in $\mathbf{b}$:
\begin{enumerate}
\item[(i)] The eight points $\lambda_{1}, \ldots, \lambda_{8}$ corresponding to the singular quadric surfaces in the family 
$Y_{\mathbf{b}} \to \mathbf{b}$.
\item[(ii)] The five points $\gamma_{1}, \ldots, \gamma_{5}$ corresponding to the singular conics in the family $\mathfrak{C} \to \mathbf{b}$.
\end{enumerate}
Now observe that each line on the quadric surface  $\mathsf{q}_{t}$ 
intersects the conic $\mathfrak{C}_{t}$ at a single point. Thus for $t \in \mathbf{b}$, the Hilbert scheme $\mycal{R}_{t}$ of lines in $\mathsf{q}_{t}$  can be described as follows:
\begin{description}
   \item[$\mycal{R}_{t} = \mycal{R}_{t^{+}}\sqcup \mycal{R}_{t^{-}} = \mathfrak{C}_{t}\sqcup \mathfrak{C}_{t}$],  when $\mathfrak{C}_{t}$ and $\mathsf{q}_{t}$ are smooth and $\pi^{-1}(t) = \{t^{+},t^{-}\} \subset C$.
   \item[$\mycal{R}_{t} = \mycal{R}_{x} = \mathfrak{C}_{t}$], when $\mathfrak{C}_{t}$ is smooth,
   $\mathsf{q}_{t}$ is singular (that is $t \in \{\lambda_{i}\}_{i = 1}^{8}$), and $x \in C$ is the ramification point sitting over the branch point $t \in \mathsf{b}$.
   \item[$\mycal{R}_{t} = \mycal{R}_{t^{+}}\sqcup \mycal{R}_{t^{-}} = \mathfrak{C}_{t}^{+}\sqcup \mathfrak{C}_{t}^{-}$], when $\mathfrak{C}_{t} = 
   \mathfrak{C}_{t}^{+}\cup \mathfrak{C}_{t}^{-}$ is reducible (that is $t \in \{\gamma_{j}\}_{j=1}^{5}$), $\mathsf{q}_{t}$ is smooth, and $\pi^{-1}(t) = \{t^{+},t^{-}\} \subset C$. 
\end{description}
This leads to the following explicit description of the Hilbert scheme $\mycal{R}|_{C}$:

Consider $\mathfrak{C}\times_{\mathbf{b}} C$. This is a smooth surface which is a conic bundle over $C$ with $10$ reducible fibers over the preimages  
$\gamma_{j}^{\pm} \in C$ of the points $\gamma_{j}$. Then $\mycal{R}|_{C}$ is obtained from $\mathfrak{C}\times_{\mathbf{b}} C$ by contracting the components $\mathfrak{C}_{\gamma_{j}}^{+}$ in the fibers over the points $\gamma_{j}^{+}$ and the components $\mathfrak{C}_{\gamma_{j}}^{-}$ in the fibers over the points $\gamma_{j}^{-}$.

On the other hand, in the surface $\mathfrak{C}$  we can  contract the components $\mathfrak{C}_{\gamma_{j}}^{+}$ in the fibers over the points $\gamma_{j} \in \mathbf{b}$. The resulting surface is a smooth $\bbP^{1}$-bundle over $\mathbf{b}$ and hence is isomorphic to $\bbP(\mathcal{W})$ for some rank two vector bundle $\mycal{W}$  on $\mathbf{b}$. By construction $\bbP(\pi^{*}\mycal{W})$ is obtained from $\mathfrak{C}\times_{\mathbf{b}} C$ by contracting the components $\mathfrak{C}_{\gamma_{j}}^{+}$  both in the fibers over $\gamma_{j}^{+}$ and in the fibers over $\gamma_{j}^{-}$. 

Therefore, the projective bundle $\mycal{R}|_{C}$ is obtained from the projective bundle $\bbP(\pi^{*}\mycal{W})$ by performing $5$ elementary modifications at points in the fibers over $\gamma_{j}^{-}$. In other words 
$\mycal{R}|_{C} \cong \bbP(W)$, where $W$ is a rank two bundle on $C$ obtained from $\pi^{*}\mycal{W}$ by $5$ simple Hecke transforms at the points $\gamma_{i}^{-}$. This implies $\deg W = \deg \pi^{*}\mycal{W} + 5 = 2\deg \mycal{W} + 5$, i.e. $\deg W$ is odd and hence by Lemma~\ref{lem:degree} we have $\fra{B\hp} = \frac{1}{2}$ as claimed.

Finally, by the argument in \cite[Lemma~6.1]{sashak-cubic} we conclude that $\fra{B^{2}}$ will be independent of the choice of $B$ and will depend only on $\alpha$. Indeed, if $u \in H^{2}(S,\bbZ)$ we will have 
\[
(B+u)^{2} = B^{2} + 2Bu + u^{2}.
\]
But $2Bu$ and $u^{2}$ are both integers and so $\fra{(B+u)^{2}} = \fra{B^{2}}$.
Similarly we have
\[
\left(B +\frac{1}{2}\hp\right)^{2} = B^{2} + B\hp + \frac{1}{2}. 
\]
Since $\fra{B\hp} = \frac{1}{2}$ we conclude that 
the fractional parts of $B^{2}$ and $\left(B +\frac{1}{2}\hp\right)^{2}$ will be equal. The corollary is proven \ \hfill $\Box$

\begin{remark}
In the proof we used the fact that a smooth conic $\mathfrak{C}$ that is a plane section of a smooth quadric surface $\mathsf{q}$ comes with a natural isomorphism to each of the two rulings $\mycal{R}^{\pm}$
on $\mathsf{q}$. In particular, it induces an isomorphism $i_\mathfrak{C}: \mycal{R}^+ \to \mycal{R}^-$
between those two rulings. In fact, the family of isomorphisms between the rulings of $\mathsf{q}$ can be identified with the family of such conics, i.e. the family of smooth curves in the linear system $\mycal{O}_{\mathsf{q}}(1,1)$: given the isomorphism $i_\mathfrak{C}$, we recover $\mathfrak{C}$ as the locus of intersections of paired lines.
\end{remark}

\

\

\section{Numerical invariants of \texorpdfstring{$S$}{S}} \label{sec:numinv}

Let $\net$ denote a copy of $\bbP^{2}$ and 
let $\pi : S \to \net$ be the double cover branched along a very general smooth curve $\discr \subset \net$  of degree $2d$. Then $S$ is a smooth simply connected Noether-Lefschetz-general surface such that $\op{Pic}(S) = \bbZ$ is generated by the hyperplane line bundle $\pi^{*}\cO_{\net}(1)$.

\subsection{Double covers of the plane} \label{ssec:doublec}
 
 \noindent
 Let $\mathsf{T} = \op{tot}(\cO_{\net}(d))$, let $p : \mathsf{T} \to \net$ be the natural projection, and let $\lambda \in H^{0}(\mathsf{T},p^{*}\cO_{\net}(d))$ be the tautological section. By definition the surface $S$ is embedded as a smooth divisor $S \subset \mathsf{T}$  in the linear system $|p^{*}\cO_{\net}(d)|$. Explicitly, if $\discr$ is the zero locus of a homogeneous polynomial  $\bdelta \in H^{0}(\net,\cO_{\net}(2d))$, then $S \subset \mathsf{T}$ is given by the equation 
\[
\lambda^{2} - p^{*}\bdelta = 0,
\]
where $\lambda^{2} - p^{*}\bdelta$ is a global section of $p^{*}\cO_{\net}(2d)$ on $\mathsf{T}$. 
In particular we have
\[
\pi_{*}\cO_{S} = \cO_{\net} \oplus \cO_{\net}(-d),
\]
and the ramification divisor $\mathsf{Ram}$ of $\pi : S \to \net$ is given by the equation 
$\lambda|_{S} = 0$. Thus $\mathsf{Ram}$ corresponds to the line bundle $p^{*}\cO_{\net}(d)|_{S} = \pi^{*}\cO_{\net}(d)$. So, by the Hurwitz formula we get
\[
\omega_{S} = \pi^{*}\omega_{\net} \otimes \cO_{S}(\mathsf{Ram}) = \pi^{*}\left( \omega_{\net}\otimes \cO_{\net}(d)\right) = \pi^{*}\cO_{\net}(d-3).
\]
By the projection formula this gives 
\[
\begin{aligned}
H^{0}(S,\omega_{S}) & = H^{0}(\net,\pi_{*}\omega_{S}) \\
& = H^{0}(\net,\cO_{\net}(d-3)\otimes \pi_{*}\cO_{S}) \\ 
& = 
H^{0}(\net,\cO_{\net}(d-3))\oplus H^{0}(\net,\cO_{\net}(-3)) \\
& = H^{0}(\net,\cO_{\net}(d-3)),
\end{aligned}
\]
i.e. 
\[
h^{0,2}(S) = h^{2,0}(S) = \binom{d-1}{2} = \frac{(d-1)(d-2)}{2}
\]

\

\subsection{Chern classes} \label{ssec:chernc}

\noindent 
Computing the Chern classes of a double cover $\pi : S \to \net$ is straightforward.  From the normal bundle sequence 
\[
0 \to T_{S} \to T_{\mathsf{T}}|_{S} \to \cO_{S}(S) \to 0
\]
of $S \subset \mathsf{T}$ we get 
\[
ch(T_{S}) = ch(T_{\mathsf{T}})\cdot S - ch(\pi^{*}\cO_{\net}(2d)).
\]
From the tangent sequence
\[
\xymatrix@R-1.5pc@C+2pc@M+0.5pc{
0 \ar[r] & T_{\mathsf{T}/\net} \ar[r] \ar@{=}[d] & T_{\mathsf{T}} \ar[r] & p^{*}T_{\net} \ar[r] & 0, \\
& p^{*}\cO_{\net}(d) & & &
}
\]
we have
\[
ch(T_{\mathsf{T}}) = p^{*}(ch(T_{\net}) + ch(\cO_{\net}(d))) = p^{*}\left(3 + (d+3) c_{1}(\cO_{\net}(1)) + \frac{d^{2} + 3}{2} \pt\right),
\]
and so 
\[
ch(T_{\mathsf{T}})\cdot S = \pi^{*}\left(3 + (d+3) c_{1}(\cO_{\net}(1)) + \frac{d^{2} + 1}{2} \mathsf{pt}\right) = 3 + (d+3) \hp + (d^{2}+ 3)\pt.
\]
Altogether we get 
\[
ch(T_{S}) = 2 + (3-d)\hp + (3-3d^{2})\pt,
\]
or
\[
c_{1}(T_{S}) = (3-d)\hp, \quad c_{2}(T_{S}) = (4d^2 -6d + 6)\pt. 
\]
In particular we have $\chi_{\op{top}}(S) = c_{2}(T_{S}) = 4d^2 -6d + 6$ and since $S$ is simply connected and has  first and third Betti numbers equal to zero, we conclude that $h^{1,1}(S) = 3d^2 - 3d + 2$.

Specializing to the case at hand when $d = 4$ we see that the map $\pi : S \to \net$ is given by the complete linear system $|\omega_{S}|$ and so  $S$ is of general type. In summary, we have proven the following

\

\begin{lemma} \label{lem:Sinvariants} Let  $\net \subset |\cO_{\bbP^{7}}(2)|$ be a very general net of quadrics, let $\pi : S \to \net$  be the double cover parametrizing the rulings in the quadrics, and let 
$h = \pi^{*}c_{1}(\cO_{\net}(1))$ denote the hyperplane class on $S$. Then:
\begin{enumerate}[label=\upshape({\bfseries \alph*}),ref=\thelemma ({\bfseries \alph*})]
\item\label{p1} $S$
is a smooth simply connected surface of general type.
\item\label{p2} $\op{Pic}(S) = \bbZ$, generated by the hyperplane line bundle $ \pi^{*}\cO_{\net}(1)$.
\item\label{p3} $\pi^{*}\cO_{\net}(1) = \omega_{S}$ and the morphism $\pi$ is given by the complete canonical linear system $|\omega_{S}|$ on $S$.
\item\label{p4} $c_{1}(S) = -\hp$, $c_{2}(S) = 46\pt$.
\item\label{p5} The Hodge diamond of $S$ is 

\

\begin{center}
\begin{tabular}{ccccc}
& & $h^{0,0}(S)$ & & \\
& $h^{0,1}(S)$ & &  $h^{1,0}(S)$ & \\
$h^{0,2}(S)$ & & $h^{1,1}(S)$ & & $h^{2,0}(S)$ \\
& $h^{1,2}(S)$ & & $h^{2,2}(S)$ & \\
& & $h^{2,2}(S)$ & & 
\end{tabular} \qquad $=$ \qquad 
\begin{tabular}{ccccc}
& & $1$ & & \\
& $0$ & &  $0$ & \\
$3$ & & $38$ & & $3$ \\
& $0$ & & $0$ & \\
& & $1$ & & 
\end{tabular}
\end{center}
\end{enumerate}
\end{lemma}

\

\section{Pairings} \label{sec:pairings}

Recall \cite{Huybrechts-FM} that if $Z$ is a smooth complex projective variety  we have a natural \emph{\bfseries Euler pairing} 
\[
\chi(-,-) \, :  \, K(Z)\otimes K(Z) \to \bbZ
\]
on the Grothendieck group  $K(Z) = K(D^{b}(Z))$, given by 
\[
\chi([\mathcal{E}],[\mathcal{F}]) = \sum_{i \in \bbZ} (-1)^{i} \dim \op{Ext}^{i}(\mathcal{E},\mathcal{F}), \qquad \text{for all} \ \mathcal{E}, \mathcal{F} \in \mathsf{ob} \, D^{b}(X).
\]
Furthermore we have the \emph{\bfseries Mukai vector} map
\[
\muk_{Z} : K(Z) \to  \oplus_{p} H^{p,p}(Z,\bbQ), \quad [\mathcal{E}] \mapsto \muk_{Z}(\mathcal{E}) = ch(\mathcal{E})\sqrt{td(Z)}
\]
which combined with the Grothendieck-Riemann-Roch statement allows us to express the Euler pairing in purely cohomological terms. Indeed, we may consider   \cite{markarian,calM1,calM2,Huybrechts-FM} the \emph{\bfseries Mukai pairing} 
\[
\langle -, -\rangle \, : \, H^{\bullet}(Z,\bbC)\otimes H^{\bullet}(Z,\bbC) \to \bbC
\]
on the  cohomology $H^{\bullet}(Z,\bbC)$, given by 
\[
\left\langle v, v'\right\rangle  \, = \, \int_{Z} (v^{\vee} \smile v')\exp(c_{1}(Z)/2),
\]
where the \emph{\bfseries dual of a cohomology class} $v = \sum_{k} v_{k} \in \oplus_{k} H^{k}(Z,\bbC)$ is
defined to be \linebreak $v^{\vee} = \sum_{k} (\iu)^{k}  v_{k}$. 

With this notation the Grothendieck-Riemann-Roch formula becomes 
\[
\chi([\mathcal{E}],[\mathcal{F}]) = \left\langle \muk_{Z}([\mathcal{E}]), \muk_{Z}([\mathcal{F}]) \right\rangle .
\]
The free abelian group $\muk_{Z}(K(Z)) \subset H^{\bullet}(Z,\bbQ)$ together with its Mukai pairing $\langle-,-\rangle$ is called the \emph{\bfseries numerical $\boldsymbol{K}$-lattice} of $Z$ and is denoted by $\Knum(Z)$.
It is isomorphic to the quotient of $K(Z)$ by the kernel of the Euler pairing and is  a derived invariant of $Z$ \cite{Huybrechts-FM}: if $Z_{1}$ and $Z_{2}$ are smooth projective varieties 
with non-isomorphic numerical $K$-lattices $\Knum(Z_{1})$ and $\Knum(Z_{2})$, then $D^{b}(Z_{1})$ and $D^{b}(Z_{2})$ can not be equivalent.

There is an obvious version $\Knum(Z,\zeta) = K(D^{b}(Z,\zeta))/\ker \chi_{(Z,\zeta)}$ of the numerical $K$-lattice which is defined for the derived category  of twisted sheaves on a smooth projective variety. When the twisting is given by a rational Brauer class $\zeta \in H^{2}(Z,\cO_{Z}^{\times})$ which is moreover topologically trivial, then the numerical $K$-lattice can again be described \cite{HStwisted} as the image in cohomology  of a twisted Mukai vector map, equipped with the Mukai pairing. We will study this invariant for the gerbe $(S,\alpha)$ corresponding to $X \subset \bbP^{7}$. 

\

\subsection{The \texorpdfstring{$\alpha$}{a}-twisted Euler pairing on \texorpdfstring{$S$}{S}} \label{ssec:alphatwpair}

\noindent
Suppose now  $\pi : S \to \net$ is a double cover branched at a very general octic and let $B \in H^{2}(S,\frac{1}{2}\bbZ)$ be a class such that $\alpha = \exp(2\pi \iu B) \in H^{2}(S,\cO_{S}^{\times})$ is a non-trivial Brauer class.

Consider the bounded derived category $D^{b}(S,\alpha)$ of $\alpha$-twisted coherent sheaves\footnote{
If we choose a sheaf $\mathcal{A} \to S$ of rank four Azumaya algebras representing $\alpha$, we can identify 
$D^{b}(S,\alpha)$ with  the derived category $D^{b}(S,\mathcal{A})$ of $\mathcal{A}$-modules. 
} on $S$. Similarly to the untwisted case we have an Euler pairing 
\[
\chi_{(S,\alpha)}(-,-)  : K(D^{b}(S,\alpha))\otimes K(D^{b}(S,\alpha)) \to \bbZ 
\]
on the Grothendieck group of $D^{b}(S,\mathcal{A})$,
given by
\[
\chi_{(S,\alpha)}([\mathcal{E}],[\mathcal{F}]) = \sum_{i} \dim \op{Ext}^{i}(\mathcal{E},\mathcal{F}), \qquad \text{for all} \ \mathcal{E}, \mathcal{F} \in \mathsf{ob} \, D^{b}(S,\alpha).
\]
By \cite{HStwisted}, the choice of $B$ as above gives a way of computing the Euler pairing $\chi_{(S,\alpha)}$ in terms of characteristic classes. More precisely:
\begin{enumerate}
\item[(1)] There exists \cite[Proposition~1.2]{HStwisted} a \emph{\bfseries $B$-twisted Chern character}
which is a linear map
\[
ch^{B} : K(D^{b}(S,\alpha)) \to H^{\bullet}(S,\bbQ).
\]
\item[(2)] The image of $ch^{B}$ is \cite[Remark~1.3(ii)]{HStwisted} a finite index subgroup in the free abelian group
\[
H(S,B;\bbZ) := \exp(B)\left(\oplus_{p} H^{p,p}(S,\bbQ)\right) \, \cap \, H^{\bullet}(S,\bbZ).
\]
\item[(3)] The $\alpha$-twisted Euler pairing is computed \cite[Remark~4.2]{HStwisted} as
\[
\chi_{(S,\alpha)}([\mathcal{E}],[\mathcal{F}]) = \left\langle \muk_{(S,B)}([\mathcal{E}]), \muk_{(S,B)}([\mathcal{F}])\right\rangle,
\]
where 
\[
\muk_{(S,B)} : \, K(D^{b}(S,\alpha)) \to H^{\bullet}(S,\bbQ), \qquad [\mathcal{E}] \mapsto ch^{B}([\mathcal{E}])\sqrt{td(S)},
\]
is the \emph{\bfseries $B$-twisted Mukai vector} map, and $\langle -, - \rangle$ is the Mukai pairing.
\end{enumerate}

\

\

\noindent
Up to isomorphism, the pair $(\muk_{(S,B)}(K(D^{b}(S,\alpha)),\langle -, - \rangle)$ depends only on $\alpha$ and only its embedding in $(H^{\bullet}(S,\bbQ),\langle-,-\rangle)$ depends on the choice of $B$ \cite{HStwisted}. 
We will write $\Knum(S,\alpha)$ for the lattice $(\muk_{(S,B)}(K(D^{b}(S,\alpha)),\langle -, - \rangle)$ and will call it the \emph{\bfseries numerical $\boldsymbol{K}$-lattice of $(S,\alpha)$}. It again 
is a derived invariant  of $(S,\alpha)$ \cite{HStwisted}.

\

\noindent
Now, following Kuznetsov's  strategy from \cite[Section~6]{sashak-cubic}, we would like to argue that 
$\Knum(S,\alpha)$ can not be isomorphic to the numerical $K$-lattice  lattice of any smooth projective surface. 
First, we have the following observation\footnote{See also the $K3$ case analyzed in \cite[Lemma~6.3]{sashak-cubic}.}:

\

\begin{lemma} \label{lem:HSB}
$H(S,B;\bbZ)$ is a finite index subgroup in the free abelian subgroup \linebreak 
$\bbZ \, (2+2B)\oplus \bbZ\, \hp \oplus \bbZ\, \pt \subset H^{\bullet}(S,\bbZ)$.
\end{lemma}

\noindent
{\bfseries Proof.} We have $H^{0}(S,\bbQ)\oplus H^{1,1}(S,\bbQ) \oplus H^{4}(S,\bbQ) = \left\{ \, r + d\hp + s\pt \, | \, r,d,s  \in \bbQ \, \right\}$, and 
\[
\exp(B)(r + d\hp + s\pt) = r + (rB + d\hp) + \left( \frac{B^{2}}{2} + d B\hp + s\right) \pt.
\]
Thus the condition $\exp(B)(r + d\hp + s\pt) \in H(S,B;\bbZ)$ is equivalent to 
\[
r \in \bbZ, \quad rB + d\hp \in H^{2}(S,\bbZ), \quad \frac{B^{2}}{2} + d B\hp + s \in \bbZ.
\]
Next note that if $r$ is odd, then $B + d\hp$ will be a class in $H^{2}(S,\bbZ)$, and so the image $\alpha$  of 
$B$ in $H^{2}(S,\cO^{\times}_{S}) = H^{2}(S,\cO_{S})/\mathsf{im} \, H^{2}(S,\bbZ)$ will be zero. But by assumption $\alpha \neq 0$ and so $r$ must be even. 

Since $r$ is even and $B \in H^{2}(S,\frac{1}{2}\bbZ)$ we will have $rB \in H^{2}(S,\bbZ)$ and so 
$d\hp \in H^{2}(S,\bbZ)$, which implies $d \in \bbZ$. Therefore 
\[
H(S,B;\bbZ) \subset \bbZ(2+2B)\oplus \bbZ\hp \oplus \bbZ\pt \subset H^{\bullet}(S,\bbZ).
\]
Since both $H(S,B;\bbZ)$ and $\bbZ(2+2B)\oplus \bbZ\hp \oplus \bbZ\pt$ are free of rank $3$, this finishes  the proof of the lemma. \ \hfill $\Box$

\

\

\subsection{The numerical \texorpdfstring{$K$}{K}-lattice of \texorpdfstring{$(S,\alpha)$}{(S,a)}} \label{ssec:MukSalpha}

\noindent
By Lemma~\ref{p4} we have $c_{1}(S) = -\hp$ and so
\[
\exp\left(\frac{c_{1}(S)}{2}\right) = \exp\left(- \frac{1}{2}\hp\right) = 1 - \frac{1}{2} \hp + \frac{1}{4}\pt.
\]
Thus if $r_{1} + \eta_{1} + d_{1}\pt,  r_{2} + \eta_{2} + d_{2}\pt \in H^{0}(S,\bbC)\oplus H^{2}(S,\bbC) \oplus H^{4}(S,\bbC) = H^{\bullet}(S,\bbC)$ are two cohomology classes on $S$, their Mukai pairing will be 
\begin{equation} \label{eq:MukaiS}
\begin{aligned}
\left\langle r_{1} + \eta_{1} + d_{1}\pt, r_{2} + \eta_{2} + d_{2}\pt \right\rangle  & = \left[ (r_{1} - \eta_{1} + d_{1}\pt)\cdot(r_{2} + \eta_{2} + d_{2}\pt)\cdot\left(1 - \frac{1}{2}\hp + \frac{1}{4}\pt\right)\right]_{4} \\[+0.5pc]
& = \frac{1}{4} r_{1}r_{2} + r_{1}d_{2} + r_{2}d_{1} - \eta_{1}\cdot\eta_{2} - \frac{1}{2}r_{1} \eta_{2}\cdot \hp + \frac{1}{2} r_{2} \eta_{1}\cdot \hp.
\end{aligned}
\end{equation}

\

\noindent 
By the previous section, the numerical $K$-lattice $\Knum(S,\alpha)$ is a finite index sublattice 
in the lattice 
\[
\left( \bbZ\, (2 + 2B) \oplus \bbZ\, \hp \oplus \bbZ\, \pt, \  \left\langle \sqrt{td(S)}\cdot(-), \sqrt{td(S)}\cdot(-)\right\rangle \right).
\]
Using \eqref{eq:MukaiS} we now compute:

\

\begin{lemma} \label{lem:Gram}
The Gram matrix of the pairing $\left\langle \sqrt{td(S)}\cdot(-), \sqrt{td(S)}\cdot(-)\right\rangle$
in the basis $2 + 2B$, $\hp$, $\pt$ is 
\begin{equation} \label{eq:Gram}
\mathbb{G} := \begin{pmatrix} 16 - 4B^{2} & -2 - 2 B\hp & 2 \\
2 - 2B\hp & -2 & 0 \\
2 & 0 & 0
\end{pmatrix}.
\end{equation}
\end{lemma}

\noindent
{\bfseries Proof.} We have
\[
\begin{aligned}
td(S) & = 1 + \frac{1}{2}c_{1}(S) + \frac{1}{12}\left( c_{1}(S)^{2} + c_{2}(S)\right)\, \pt \\[+0.7pc]
& = 1 - \frac{1}{2}\hp + \frac{1}{12}\left( (-\hp)^{2}
 + 46\right)\, \pt \\[+0.7pc]
 & = 1 - \frac{1}{2} \hp + 4\, \pt.
 \end{aligned}
\]
Therefore
\[
\sqrt{td(S)} = 1 - \frac{1}{4}\hp + \frac{31}{16}\pt,
\]
and so
\begin{equation} \label{eq:normalizedbasis}
\begin{aligned}
\sqrt{td(S)}(2+2B) & = 2 + \left(2B - \frac{1}{2} \hp\right) + \left( \frac{31}{8} - \frac{1}{2}B\hp \right)\, \pt, \\[+1pc]
\sqrt{td(S)}\hp & =  0 + \hp - \frac{1}{2}\, \pt, \\[+1pc]
\sqrt{td(S)}\pt & = 0 + 0 + \pt.
\end{aligned}
\end{equation}
Plugging these into the formula \eqref{eq:MukaiS} yields the Gram matrix \eqref{eq:Gram}. \ \hfill $\Box$

\

\

\noindent
Next we compute the fractional part $\fra{B^{2}}$:

\

\begin{lemma} \label{lem:{B^2}}
$\fra{B^{2}} = \frac{3}{4}$.
\end{lemma}

\noindent
{\bfseries Proof.} Let $\mathcal{A}$ be the sheaf of  Azumaya algebras representing the gerbe $(S,\alpha)$. Then on $S$ we have an $\alpha$-twisted  
rank $2$ vector bundle $\mathsf{W}$ so that $\mycal{R} = \bbP(\mathsf{W})$ and $\mathcal{A} = \mathsf{W}^{\vee}\otimes \mathsf{W}$, where $\mathsf{W}^{\vee}$ denotes the dual $\alpha^{-1}$-twisted rank two bundle. Furthermore by \cite{sashak-quadrics} we have  an isomorphism $\pi_{*}\mathcal{A} \cong \mathcal{B}_{0}$ of sheaves of algebras on $\net$, where $\mathcal{B}_{0}$ is the even part of the sheaf of Clifford algebras for the family of quadrics $\sigma : E \to E^{\vee}\otimes \cO_{\net}(1)$ and moreover we have an equivalence $D^{b}(S,\alpha) = D^{b}(S,\mathcal{A}) \cong D^{b}(\net,\mathcal{B}_{0})$ between the derived category of $\mathcal{A}$-modules on $S$ 
and  the derived category of $\mathcal{B}^{0}$-modules on $\net$.

Now by the discussion of twisted Chern classes above we have 
\[
ch^{B}(\mathsf{W}) = 2 + (2B + d\hp) + s\pt
\]
for some integers $d, s \in \bbZ$. (Note: these are not the same as $d,s$ in Lemma \ref{lem:HSB}, where $s$ was rational.) Also by the basic properties \cite{HStwisted} of $ch^{B}$ we will have  
\[
ch^{-B}(\mathsf{W}^{\vee}) = 2 - (2B + d\hp) + s\pt,
\]
and
\[
\begin{aligned}
ch(\mathcal{A}) & = ch^{-B}(\mathsf{W}^{\vee})ch^{B}(\mathsf{W}) \\
& = (2 + (2B + d\hp) + s\pt)\cdot (2 - (2B + d\hp) + s\pt) \\
& = 4 + 0 + (4s - (2B+d\hp)^{2})\pt.
\end{aligned}
\]
Write $a = 4s - (2B+d\hp)^{2} \in \bbZ$. Then
since we computed 
\[
td(S) = 1 - \frac{1}{2}\, \hp + 4\, \pt,
\]
we get 
\[
ch(\mathcal{A})td(S) = 4 - 2\, \hp + (16 + a)\, \pt,
\]
and hence
\[
\pi_{*}\left(ch(\mathcal{A})td(S)\right) = 8 - 4\, \Hp  + (16 + a)\, \pt,
\]
where $\Hp \in H^{2}(\net,\bbZ)$ is the hyperplane class on $\net$.

Also, as a sheaf we have 
\[
\pi_{*}\mathcal{A} = \mathcal{B}_{0} = \cO_{\net} \oplus \wedge^{2}E\otimes \cO_{\net}(-1) \oplus \wedge^{4}E \otimes \cO_{\net}(-2).
\]
From the short exact sequence 
\[
0 \to E \to \cO_{\net}^{\oplus 6} \to \cO_{\net}(1)^{\oplus 2} \to 0
\]
we compute
\[
\wedge^{4} E = \cO_{\net}(-2), \qquad c_{1}(E) = - 2\, \Hp, \ ch_{2}(E) = -\pt.
\]
Since 
\[
ch(\wedge^{2}E) = 6 + 3 c_{1}(E) + \left( 2ch_{2}(E) + \frac{c_{1}(E)^{2}}{2}\right)\, \pt,
\]
we get 
\[
ch(\wedge^{2}E) = 6 - 6\, \Hp + 0\, \pt,
\]
and hence
\[
\begin{aligned}
ch(\wedge^{2}E\otimes \cO_{\net}(-1)) & = (6 - 6\, \Hp + 0\, \pt)\cdot 
\left(1 - \Hp + \frac{1}{2}\, \pt \right) \\
& = 6 - 12\, \Hp + 9\, \pt. 
\end{aligned}
\]
Similarly
\[
ch(\cO_{A}(-4)) = 1 - 4\, \Hp + 8\, \pt,
\]
and so 
\[
\begin{aligned}
ch(\pi_{*}\mathcal{A}) & = ch(\cO_{\net}) + ch(\wedge^{2}E\otimes \cO_{\net}(-1)) + ch(\cO_{A}(-4))  \\
& = 8 - 16\, \Hp + 17\, \pt.
\end{aligned}
\]
Since $\net \cong \bbP^{2}$ we have 
\[
td(\net) = 1 + \frac{3}{2}\, \Hp + \pt,
\]
and thus
\[
ch(\pi_{*}\mathcal{A})td(\net) = 8 - 4\, \Hp + pt.
\]
Therefore by  Grothendieck-Riemann-Roch we get 
\[
16 + a = 1 
\]
or 
\[
4s - (2B + d\, \hp)^{2}  + 16 = 1. 
\]
But
\[
(2B + d\, \hp)^{2} = 4B^{2} + 4d\, B\hp + 2d^{2} = 4B^{2} + 2d(d + 2 B\hp).
\]
But $2B\hp$ is odd, by Corollary \ref{cor:{Bh}}, and so $2d(d+B\hp)$ is divisible by $4$. In other words
$4B^{2} = -1$ modulo $4$, and so $\fra{B^{2}} = \frac{3}{4}$.  \ \hfill $\Box$

\

\

\begin{corollary} \label{cor:obstruction}
Suppose that $\fra{B\hp} = \frac{1}{2}$ and $\fra{B^{2}} = \frac{3}{4}$. Then there are no elements  $\mathsf{v}_{1}, \mathsf{v}_{2} \in \Knum(S,\alpha)$ with $\chi_{(S,\alpha)}(\mathsf{v}_{1}, \mathsf{v}_{2}) = 1$ and $\chi_{(S,\alpha)}(\mathsf{v}_{2}, \mathsf{v}_{2}) = 0$.
\end{corollary}

\noindent
{\bfseries Proof.}  By Lemma~\ref{lem:HSB} and Lemma~\ref{lem:Gram} the numerical $K$-lattice $\Knum(S,\alpha)$ is a finite index sublattice in the lattice $\bbZ^{\oplus 3}$ with a pairing given by the Gram matrix \eqref{eq:Gram}.
Thus it suffices to show that there are no elements $\mathsf{v}_{1}, \mathsf{v}_{2} \in \bbZ^{\oplus 3}$ for which 
\[
\mathsf{v}_{1}^{T} \mathbb{G} \mathsf{v}_{2} = 1, \qquad \text{and} \ \mathsf{v}_{2}^{T} \mathbb{G} \mathsf{v}_{2} = 0.
\]

Let $\mathsf{v}_{2} = \begin{pmatrix} x \\ y \\ z \end{pmatrix} \in \bbZ^{\oplus 3}$. The condition 
$\mathsf{v}_{2}^{T} \mathbb{G} \mathsf{v}_{2} = 0$ reads
\[
(16 - 4B^{2})x^{2} - 2 y^{2}  - 4B\hp xy + 4 xz = 0.
\]
Reducing mod 4, this becomes
\[
x^{2} - 2 y^{2}  - 2 xy  = 0 \ \text{mod} \ 4.
\]
Here all terms are even except for the first term, so $x$ must be even. But then all terms are divisible by 4 except for the second term, so $y$ too must be even.
But then $ \mathbb{G} \mathsf{v}_{2} $ is even. This contradicts the requirement that $\mathsf{v}_{1}^{T} \mathbb{G} \mathsf{v}_{2} = 1$

This completes the proof of  the corollary. \ \hfill $\Box$

\

\noindent
Putting all of this together   we get the following

\ 

\begin{theorem} \label{theo:nc}
Let $X$ be a very general intersection of $3$-quadrics in $\bbP^{7}$, and let 
$(S,\alpha)$ be the associated surface of general type with a $2$-torsion Brauer class. Then the associated numerical $K$-lattice $(\Knum(S,\alpha),\chi_{(S,\alpha)})$ is not isomorphic to $(\Knum(M),\chi_{M})$ for any smooth projective  surface $M$. In particular $D^{b}(S,\alpha)$ can not be equivalent to the derived category of any smooth projective surface.  
\end{theorem}

\noindent
{\bfseries Proof:} \ If $M$ is a smooth projective surface and $x \in M$, then 
we have $\chi_{M}([\cO_{M}],[\cO_{x}]) = 1$, and $\chi_{M}([\cO_{x}],[\cO_{x}]) = 0$. But by Lemma~\ref{lem:{B^2}}, Corollary~\ref{cor:{Bh}}, and 
Corollary~\ref{cor:obstruction} there are no elements $\mathsf{v}_{1}, \mathsf{v}_{2} \in \Knum(S,\alpha)$, satisfying $\chi_{(S,\alpha)}(\mathsf{v}_{1},\mathsf{v}_{2}) = 1$, and $\chi_{(S,\alpha)}(\mathsf{v}_{2},\mathsf{v}_{2}) = 2$. Thus $(\Knum(S,\alpha),\chi_{(S,\alpha)}) \ncong (\Knum(M),\chi_{M})$, and $D^{b}(S,\alpha) \ncong D^{b}(M)$. \ \hfill $\Box$

\

\bibliographystyle{halpha}
\bibliography{quadrics}

\end{document}